# A Generalized Numeration Base

by Florentin Smarandache, Ph. D.
University of New Mexico
Gallup, NM 87301, USA

**Abstract.** Á Generalized Numeration Base is defined in this paper, and then particular cases are presented, such as Prime Base, Square Base, m-Power Base, Factorial Base, and operations in these bases.

## Introduction.

The following bases are important for partitions of integers into primes, squares, cubes, generally into m-powers, also into factorials, and into any strictly increasing sequence.

### 1) The Prime Base:

0,1,10,100,101,1000,1001,10000,10001,10010,10100,100000,100001,1000000, 1000001,1000010,1000100,10000000,10000001,100000000,100000001,100000010, 100000100,1000000000,1000000001,1000000010,1000000100,1000000101,... .
(Each number n written in the Prime Base.)

(We define over the set of natural numbers the following infinite base: $p_0 = 1$, and for $k \geq 1$ $p_k$ is the k-th prime number.)
He proved that every positive integer A may be uniquely written in the Prime Base as:

$$A = \overline{(a_n \ldots a_1 a_0)}_{(SP)} \stackrel{def}{===} \sum_{i=0}^{n} a_i p_i, \text{ with all } a_i = 0 \text{ or } 1, \text{ (of course } a_n = 1\text{),}$$

in the following way:
- if $p_n \leq A < p_{n+1}$ then $A = p_n + r_1$;
- if $p_m \leq r_1 < p_{m+1}$ then $r_1 = p_m + r_2$, $m < n$;
and so on until one obtains a rest $r_j = 0$.

Therefore, any number may be written as a sum of prime numbers + e, where e = 0 or 1.

If we note by p(A) the superior prime part of A (i.e. the largest
prime less than or equal to A), then
A is written in the Prime Base as:

   A = p(A) + p(A-p(A)) + p(A-p(A)-p(A-p(A))) + ... .

This base is important for partitions with primes.

## 2) The Square Base:
0,1,2,3,10,11,12,13,20,100,101,102,103,110,111,112,1000,1001,1002,1003,
1010,1011,1012,1013,1020,10000,10001,10002,10003,10010,10011,10012,10013,
10020,10100,10101,100000,100001,100002,100003,100010,100011,100012,100013,
100020,100100,100101,100102,100103,100110,100111,100112,101000,101001,
101002,101003,101010,101011,101012,101013,101020,101100,101101,101102,
1000000,... .
(Each number n written in the Square Base.)

(We define over the set of natural numbers the following infinite
base:  for k >= 0   $s_k = k^2$.)

We prove that every positive integer A may be uniquely written in
the Square Base as:

$$A = \overline{(a_n \ldots a_1 a_0)}_{(S2)} \overset{def}{===} \sum_{i=0}^{n} a_i s_i, \text{ with } a_i = 0 \text{ or } 1 \text{ for } i \geq 2,$$

$0 \leq a_0 \leq 3$, $0 \leq a_1 \leq 2$, and of course $a_n = 1$,
in the following way:
  - if $s_n \leq A < s_{n+1}$  then $A = s_n + r_1$ ;
  - if $s_m \leq r_1 < p_{m+1}$ then $r_1 = s_m + r_2$, m < n;
  and so on until one obtains a rest $r_j = 0$.

Therefore, any number may be written as a sum of squares (1 not counted
as a square -- being obvious) + e, where e = 0, 1, or 3.

If we note by s(A) the superior square part of A (i.e. the
largest square less than or equal to A), then A is written in the
Square Base as:

   A = s(A) + s(A-s(A)) + s(A-s(A)-s(A-s(A))) + ... .

This base is important for partitions with squares.

## 3) The m-Power Base (generalization):

(Each number n written in the m-Power Base,
where m is an integer >= 2.)

(We define over the set of natural numbers the following infinite
m-Power Base: for k >= 0  $t_k$ = k^m.)

He proved that every positive integer A may be uniquely written in
the m-Power Base as:

$$A = \overline{(a_n \ldots a_1 a_0)}_{(SM)} \stackrel{def}{===} \sum_{i=0}^{n} a_i t_i, \text{ with } a_i = 0 \text{ or } 1 \text{ for } i \geq m,$$

$0 \leq a_i \leq \lfloor ((i+2)^m - 1) / (i+1)^m \rfloor$ (integer part)

for i = 0, 1, ..., m-1, $a_i$ = 0 or 1 for i >= m, and of course $a_n$ = 1,
in the following way:
- if $t_n \leq A < t_{n+1}$  then $A = t_n + r_1$ ;
- if $t_m \leq r_1 < t_{m+1}$  then $r_1 = t_m + r_2$ , m < n;
and so on until one obtains a rest $r_j$ = 0.

Therefore, any number may be written as a sum of m-powers (1 not counted
as an m-power -- being obvious) + e, where e = 0, 1, 2, ..., or 2^m-1.

If we note by t(A) the superior m-power part of A (i.e. the
largest m-power less than or equal to A), then A is written in the
m-Power Base as:

   A = t(A) + t(A-t(A)) + t(A-t(A)-t(A-t(A))) + ...

This base is important for partitions with m-powers.

## 4) The Factorial Base:
0,1,10,11,20,21,100,101,110,111,120,121,200,201,210,211,220,221,300,301,310,
311,320,321,1000,1001,1010,1011,1020,1021,1100,1101,1110,1111,1120,1121,
1200,...
(Each number n written in the Factorial Base.)

(We define over the set of natural numbers the following infinite
base: for k >= 1  $f_k$ = k!)

He proved that every positive integer A may be uniquely written in
the Factorial Base as:

$$A = (a_n \ldots a_2 a_1)_{(F)} \equiv \sum_{i=1}^{} a_i f_i, \text{ with all } a_i = 0, 1, \ldots, i \text{ for } i \geq 1.$$

in the following way:
- if $f_n \leq A < f_{n+1}$ then $A = f_n + r_1$;
- if $f_m \leq r_1 < f_{m+1}$ then $r_1 = f_m + r_2$, $m < n$;

and so on until one obtains a rest $r_j = 0$.

What's very interesting: $a_1 = 0$ or 1; $a_2 = 0, 1,$ or 2; $a_3 = 0, 1, 2,$ or 3, and so on...

If we note by f(A) the superior factorial part of A (i.e. the largest factorial less than or equal to A), then A is written in the Factorial Base as:

   A = f(A) + f(A-f(A)) + f(A-f(A)-f(A-f(A))) + ... .

Rules of addition and subtraction in Factorial Base:
For each digit $a_I$ we add and subtract in base i+1, for i >= 1.

For example, an addition:
```
                    base  5 4 3 2
                         ---------------
                            2 1 0 +
                            2 2 1
                         -----------
                          1 1 0 1
```
because: 0+1= 1   (in base 2);
         1+2=10   (in base 3), therefore we write 0 and keep 1;
         2+2+1=11 (in base 4).

Now a subtraction:
```
                    base  5 4 3 2
                         ---------------
                            1 0 0 1 -
                              3 2 0
                         ---------
                          = = 1 1
```
because: 1-0=1 (in base 2);
         0-2=? it's not possible (in base 3),
             go to the next left unit, which is 0 again (in base 4),
             go again to the next left unit, which is 1 (in base 5),
             therefore 1001 --> 0401 --> 0331
             and then 0331-320=11.

Find some rules for multiplication and division.

In a general case:
   if we want to design a base such that any number

$$A = (a_n \ldots a_2 a_1)_{(B)} \stackrel{def}{===} \sum_{i=1}^{n} a_i b_i \text{, with all } a_i = 0, 1, \ldots, t_i \text{ for}$$

$i \geq 1$, where all $t_i \geq 1$, then:
this base should be
$b_1 = 1$, $b_{i+1} = (t_i + 1) * b_i$ for $i \geq 1$.

## 5) The Generalized Numeration Base:

(Each number n written in the Generalized Numeration Base.)

(We define over the set of natural numbers the following infinite Generalized Numeration Base: $1 = g_0 < g_1 < \ldots < g_k < \ldots$ .)

He proved that every positive integer A may be uniquely written in the Generalized Numeration Base as:

$$A = (a_n \ldots a_1 a_0)_{(SG)} \stackrel{def}{===} \sum_{i=0}^{n} a_i g_i \text{, with } 0 \leq a_i \leq \lfloor (g_{i+1} - 1) / g_i \rfloor$$

(integer part) for $i = 0, 1, \ldots, n$, and of course $a_n \geq 1$,
in the following way:
- if $g_n \leq A < g_{n+1}$ then $A = g_n + r_1$;
- if $g_m \leq r_1 < g_{m+1}$ then $r_1 = g_m + r_2$, $m < n$;
and so on until one obtains a rest $r_j = 0$.

If we note by g(A) the superior generalized part of A (i.e. the largest $g_i$ less than or equal to A), then A is written in the Generalized Numeration Base as:

$A = g(A) + g(A-g(A)) + g(A-g(A)-g(A-g(A))) + \ldots$

This base is important for partitions: the generalized base may be any infinite integer set (primes, squares, cubes, any m-powers, Fibonacci/Lucas numbers, Bernoully numbers, Smarandache sequences, etc.) those partitions are studied.

A particular case is when the base verifies: $2g_i \geq g_{i+1}$ for any i,

and $g_0 = 1$, because all coefficients of a written number in this base will be 0 or 1.

Remark: another particular case: if one takes $g_i = p^{i-1}$, $i = 1, 2, 3, \ldots$, p an integer $\geq 2$, one gets the representation of a number in the numerical base p {p may be 10 (decimal), 2 (binary), 16 (hexadecimal), etc.}.